\documentclass[11pt]{article}
\usepackage[leqno]{amsmath} 
\usepackage{amssymb,amscd,amsthm}
\input xy \xyoption{all}
\CompileMatrices
\def\eps{{\epsilon}}
\def\im{\operatorname{im}}
\def\Cl{\operatorname{Cl}}

\def\Oh{{\cal O}}
\def\Gh{{\cal G}}
\def\Hh{{\cal H}}

\def\LG{\Lambda (\Gh) }
\def\LH{\Lambda (\Hh) }

\def\LGamma{\Lambda (\Gamma)}
\def\N{\mathbb{N}}
\def\Z{\mathbb{Z}}
\def\Q{\mathbb{Q}}
\def\C{\mathbb{C}}
\def\R{\mathbb{R}}
\newcommand{\La}{\mathbb{L}}

\def\Kbar{{\bar{K}}}
\def\Xbar{{\bar{X}}}
\def\Kcyc{K^{\operatorname{cyc}}}

\def\Lcyc{L^{\operatorname{cyc}}}

\def\isom{\cong}
\def\prolim{\varprojlim}
\def\Aut{\operatorname{Aut}}
\def\Hom{\operatorname{Hom}}
\def\Gal{\operatorname{Gal}}
\def\Gl{\operatorname{Gl}}
\def\Tor{\operatorname{Tor}}
\def\Tw{\operatorname{Tw}}
\def\tr{\operatorname{tr}}
\def\pr{\operatorname{pr}}

\def\id{\operatorname{id}}

\def\cores{\operatorname{cores}}

\def\id{\operatorname{id}}
\def\dim{\operatorname{dim}}

\def\Spec{\operatorname{Spec}}

\def\qed{\hfill $\Box$ \medskip \par}

\newcommand{\bew}{\begin{proof}}
\newcommand{\bewende}{\end{proof}}
\newtheorem{lemma}{Lemma}[subsection]
\newtheorem{prop}[lemma]{Proposition}
\newtheorem{thm}[lemma]{Theorem}
\newtheorem{defn}[lemma]{Definition}
\newtheorem{cor}[lemma]{Corollary}

\newcommand{\rem}{\noindent{\bf Remark:\ }}

\newtheorem{conj}[lemma]{Conjecture}

\title{\vskip -25mm On non-commutative twisting in {\'e}tale and
motivic cohomology}   
\author{Jens Hornbostel and Guido Kings}
\date{April 13, 2004}

\begin{document}

\maketitle

\begin{abstract}It is proved that under a technical condition the
{\'e}tale cohomology groups $H^1(\Oh_K[1/S],H^i(\Xbar,\Q_p(j))$,
where $X\to \Spec\Oh_K[1/S]$ is a smooth, projective scheme, are generated
by twists of norm compatible units in a tower of number fields associated to
$H^i(\Xbar,\Z_p(j))$. 
This confirms a consequence of the non-abelian Iwasawa main
conjecture. Using the ``Bloch-Kato-conjecture'' a similar result is
proven for motivic cohomology with finite coefficients.
\end{abstract}

\section*{Introduction} One of the most astonishing consequences of
the equivariant Tamagawa number conjecture is the twist invariance 
of the zeta elements, which implies that 
all motivic elements should be twists
of norm compatible units in (big) towers of number fields.
More precisely
one expects that for a $\Z_p$-lattice $T$ in a motive with $p$-adic realization
$V$ the image of the twisting map (see \ref{twistingmap} below) 
\[
(\prolim_n H^1(\Oh_{K_n}[1/S], \Z_p(1)))\otimes T(j-1)
\to  H^1(\Oh_{K}[1/S], T(j))
\]
generates a subgroup of finite index. Here the inverse limit
runs over the number fields $K_n:=K(T/p^n)$ obtained
from $K$ by adjoining the elements
$T/p^n$. Moreover, the image of this map
should have a motivic meaning, that is the elements should be in
the image of the $p$-adic regulator from motivic cohomology. 
(This is explained in \cite{Hu-Ki2} and builds
on ideas of Kato \cite{Kato}).

The philosophy of 
twisting originates from work of Iwasawa, Tate and Soul{\'e}, who
considered twisting with the cyclotomic character. This
already lead to many interesting results. Here Kato's work \cite{Kato2}
on the Birch-Swinnerton-Dyer conjecture is the most 
spectacular example. Earlier Soul{\'e} used this idea in the case 
of Tate motives in his investigations about the connection
of $K$-theory and {\'e}tale cohomology for number rings \cite{So0}. 
He also pointed the way
to applications to CM-elliptic curves \cite{So2}. 

The goal of this paper is to show that the above twisting map
has for $j>>0$ indeed finite cokernel assuming
the very reasonable condition that the Iwasawa $\mu$-invariant of the
number field $K$ vanishes.
In the second part of the paper we consider the statement, that
the resulting elements are in the image of the regulator from
motivic cohomology. 
Our results in this direction give a weak hint that the elements obtained
as twists of units are motivic.
Using the ``Bloch-Kato-conjecture'' for all fields (as announced by 
Voevodsky), we prove that
there is a twisting map for motivic cohomology compatible with the
one for {\'e}tale cohomology under the cycle class map.\\

The authors like to thank J. Coates and R. Sujatha for useful discussions
and for making available the results of \cite{C-S} before their 
publication. The authors are indebted to T. Geisser for insisting
to use motivic cohomology with finite coefficients instead of 
$K$-theory in the formulation of the results in the second part.

\section{Non-commutative twisting in {\'e}tale cohomology}
In this section we describe the {\'e}tale situation.
All cohomology groups in this paper are  
{\'e}tale cohomology groups unless explicitly labeled
otherwise.

\subsection{The twisting map in {\'e}tale cohomology} \label{sectiontwist}
Let $K$ be a number field with ring of integers $\Oh_K$. Fix a
prime number $p>2$ and a finite set of primes $S$ of $\Oh_K$, which
contains the primes dividing $p$. As usual let $G_{S}:=\Gal(K_S/K)$ be the
Galois group of $K_S/K$, where $K_S$ is
the maximal outside of $S$ unramified extension field in a fixed 
algebraic closure $\Kbar$ of $K$. Let $T$
be a finitely generated free $\Z_p$-module with a continuous $G_{S}$-action
\[
\rho: G_S \longrightarrow \Aut_{\Z_p}(T).
\]
We will consider $T$ also as {\'e}tale sheaf on $\Oh_K[1/S]$ 
(see e.g. \cite[p. 640]{FPR} using the fixed embedding into $\Kbar$
as base point) and write $V:=T\otimes_{\Z_p}\Q_p$ for the 
associated $G_S$-representation and $\Q_p$-sheaf.
Let $\Gh:=\im \rho $ be the image of $\rho$. If we define finite groups
\[
G_n:= \im \left\{ \rho_n: G_S\longrightarrow  \Aut_{\Z/ p^n\Z}(T/p^nT)\right\}
\]
then we also have $\Gh\isom \prolim_n G_n$. Note that $\Gh$ is a $p$-adic
Lie group. The {\em Iwasawa algebra} of $\Gh$ is by definition
the continuous group ring
\[
\LG:=\prolim_n \Z_p[G_n] \cong \prolim_n \Z/p^n[G_n].
\]
We denote by $K_\infty$ the field fixed by the kernel of $\rho$
and by $K_n$ the field fixed by the kernel of $\rho_n$, so that
$K_\infty = \cup_n K_n$ and $\Gal(K_\infty / K)\isom \Gh$.
Note that $\Oh_{K_n}[1/S]$ is finite and etale over $\Oh_K[1/S]$.\\

\noindent{\bf Example:} To make the above definitions more concrete, consider
the following important example. Let $E/K$ be an elliptic curve
without complex multiplication and $T_pE:=\prolim_nE[p^n]$ its Tate-module. 
We have $K_n:=K(E[p^n])$.
It is a well-known result of Serre that the image of 
the Galois group $G_S$ in $\Aut_{\Z_p}(T_pE)$ has finite index and is
equal to $\Aut_{\Z_p}(T_pE)$ for almost all $p$. If we assume the latter case,
we have in the above notation $G_n\isom\Gl_2(\Z/p^n)$ and 
$\Gh\isom\Gl_2(\Z_p)$. 

The following proposition for the {\'e}tale cohomology should be well-known. For 
convenience of the reader and to explain the normalizations of the
action in detail, we give the proof in an appendix.
\begin{prop}\label{twistisom}(see appendix \ref{twistbew})
There are canonical isomorphisms of compact finitely generated $\LG$-modules
\[
H^i(\Oh_K[1/S], \LG)\otimes_{\Z_p}T\xrightarrow{\cong}H^i(\Oh_K[1/S], \LG\otimes_{\Z_p}T)
\]
and 
\[ \prolim_nH^i(\Oh_{K_n}[1/S],T/p^n) \isom
\prolim_nH^i(\Oh_{K_n}[1/S],T)\isom H^i(\Oh_K[1/S], \LG\otimes_{\Z_p}T)
\]
(limit over the corestriction maps).
Here the $\LG$-module structure on $T$ is induced by the action of $\Gh$ and the one
on $\LG$ is via multiplication
with the inverse.
\end{prop}

\rem If one does not want to assume that $T$ is a free $\Z_p$-module,
then one should use the derived tensor product in the formulas above.
This is helpful, if one wants to consider $T:=H^i(\Xbar,Z_p)$
for a smooth, projective variety $X$ over $\Oh_K[1/S]$.

Denote by $\epsilon :\LG\to \Z_p$ the augmentation map. 
\begin{defn}\label{twistingmap} The twisting map 
\[
\Tw_T:H^i(\Oh_K[1/S], \LG)\otimes_{\Z_p}T\longrightarrow H^i(\Oh_K[1/S], T)
\]
is the composition of the isomorphism of proposition
\ref{twistisom} with the map 
\[
\epsilon:H^i(\Oh_K[1/S], \LG\otimes_{\Z_p}T)\to H^i(\Oh_K[1/S], T)
\]
induced by the augmentation $\epsilon$.
\end{defn}
Our goal is to show that the twisting map is surjective in certain
cases after
tensoring with $\Q_p$. In particular it allows to construct elements in 
$H^i(\Oh_K[1/S], T)$ starting from (corestriction or norm compatible) elements in 
\[
H^i(\Oh_K[1/S], \LG)\isom \prolim_n H^i(\Oh_{K_n}[1/S],\Z_p).
\]
We will apply this in the case where 
$T=H^r(X\times_{\Oh_K[1/S]} \Kbar, \Z_p)$. 

\subsection{The conjectured ranks of the {\'e}tale cohomology}
For the convenience of the reader, we recall the conjecture
of Jannsen \cite{Ja} about the ranks of the {\'e}tale cohomology.

Let $X$ be a smooth, projective scheme over $\Oh_K[1/S]$
and denote by $\Xbar:=X\times_{\Oh_K[1/S]}\Kbar$ the base change to 
the algebraic closure. 
\begin{conj}(Jannsen)\label{Jannsen} For $i+1<j$ or $i+1>2j$ one has
\[
H^2(\Oh_K[1/S], H^i(\Xbar,\Q_p(j)))=0.
\]
\end{conj}
As a consequence one obtains (for $p\neq 2$) 
the following formula for the dimension
of the $H^1$: for $i+1<j$
\[
\dim_{\Q_p}H^1(\Oh_K[1/S], H^i(\Xbar,\Q_p(j)))=
\dim_\R H^i(X\times_\Q\C,\R(j))^+,
\]
where ``$+$'' denotes the invariants under complex conjugation,
which acts on $X\times_\Q\C$ and on $\R(j)=(2\pi i)^j\R$.

Moreover in analogy with Beilinson's conjecture that
the regulator from $K$-theory to Beilinson-Deligne cohomology is an
isomorphism for $i+1<j$, Jannsen also conjectures that the
Soul{\'e} regulator 
\[
r_p:H^{i+1}_{mot}(X,\Z(j))\otimes_\Z\Q_p\to 
H^1(\Oh_K[1/S], H^i(\Xbar,\Q_p(j)))
\]
is an isomorphism for $i+1<j$. It is shown in \cite{So2} that for 
$p$-adic $K$-theory the above regulator is surjective if $j>>0$.
This should be compared with the result in \ref{twistcor}.

\subsection{Another description of the  twisting map}\label{twistvariant}
To make the similarity with the twisting map in motivic cohomology 
and in $p$-adic $K$-theory
more apparent, we describe  the twisting map at finite level.

Fix an integer $n>0$. We have by Shapiro's lemma
\[
H^i(\Oh_K[1/S],\Z/p^n\Z[G_n])\isom H^i(\Oh_{K_n}[1/S],\Z/p^n\Z)
\]
using the identification as explained in appendix \ref{twistbew}. 
As $T/p^nT$ is a trivial sheaf over
$\Oh_{K_n}[1/S]$, we have $T/p^nT\isom H^0(\Oh_{K_n}[1/S],T/p^nT)$
and the cup product gives an isomorphism
\[
H^i(\Oh_{K_n}[1/S],\Z/p^n\Z)\otimes T/p^nT\xrightarrow{\cup}
H^i(\Oh_{K_n}[1/S],T/p^nT).
\]
Together with the corestriction (=trace map in {\'e}tale cohomology) 
\[
H^i(\Oh_{K_n}[1/S],T/p^nT)\to H^i(\Oh_{K}[1/S],T/p^nT)
\]
we get a map
\begin{equation}\label{twistatfinitelevel}
H^i(\Oh_K[1/S],\Z/p^n\Z[G_n])\otimes T/p^nT\to 
H^i(\Oh_{K}[1/S],T/p^nT).
\end{equation}

Observe that by Mittag-Leffler we have 
$\prolim_n H^i(\Oh_{K}[1/S],T/p^nT) \cong H^i(\Oh_{K}[1/S],T)$.

\begin{lemma}\label{twistiscup}
The inverse limit with respect to the trace map and
reduction on the coefficients of the maps (\ref{twistatfinitelevel})
coincides with the twisting map in definition \ref{twistingmap}.
\end{lemma}
\bew Straightforward.
\bewende

\subsection{Tate twist}
Let $\Kcyc:=\cup_nK(\mu_{p^n})$ be the field $K$ with all 
the $p$-th power roots
of unity $\mu_{p^\infty }$ adjoined. We will assume that
$K_\infty$ contains $\Kcyc$. If this is not the case it can be achieved
by considering $K_n(\mu_{p^n})$ instead of $K_n$. Let $\Gamma :=\Gal(\Kcyc /K)$, then
we have a map $\Gh\to \Gamma$ and we denote its kernel by $\Hh$. This map induces
a surjection
\begin{equation}\label{groupsurj}
\LG\to \LGamma.
\end{equation}
The cyclotomic character induces an inclusion of
$\Gamma$ in $\Z_p^*$ and the associated free $\Z_p-\Gamma$-module of 
rank $1$ is denoted by $\Z_p(1)$.  
As usual let $\Z_p(j):=\Z_p(1)^{\otimes j}$ and
$T(j):=T\otimes \Z_p(j)$.

We will consider the following important variant of the twisting map,
given by combining Definition \ref{twistingmap} and Proposition 
\ref{twistisom} for $T=\Z_p(1)$ and $T(j-1)$:
\begin{equation}\label{twistingvariant}
H^i(\Oh_K[1/S], \LG(1))\otimes_{\Z_p}T(j-1)\xrightarrow{\Tw_{T(j-1)}} 
H^i(\Oh_K[1/S], T(j)).
\end{equation}
Note that $H^1(\Oh_K[1/S], \LG(1))=\prolim_n H^1(\Oh_{K_n}[1/S], \Z_p(1))$
and that by Kummer theory we have an exact sequence
\begin{equation}
0\to \Oh_{K_n}[1/S]^*\otimes\Z/p^n\to  
H^1(\Oh_{K_n}[1/S], \Z/p^n(1))\to \Cl(\Oh_{K_n}[1/S])[p^n]\to 0.
\end{equation}
Here $\Cl(\Oh_{K_n}[1/S])[p^n]$ is the $p^n$-torsion subgroup of
the class group of $\Oh_{K_n}[1/S]$. Taking the limit over $n$ we
get 
\begin{equation}\label{Kummer}
\prolim_n \Oh_{K_n}[1/S]^*\otimes\Z_p\isom\prolim_n 
H^1(\Oh_{K_n}[1/S], \Z_p(1)).
\end{equation}

We have a twisted variant of Lemma \ref{twistiscup}.
Namely, the map $Tw_{T(j-1)}$ of (\ref{twistingvariant})
is again given by taking cup products 
$H^i(\Oh_{K_n}[1/S],\Z/p^n(1))\otimes_{\Z_p}
T/p^n(j-1)$ and passing to the limit, using again Proposition \ref{twistisom}.

\subsection{The cokernel of the twisting map}
To study the cokernel of the twisting map, we 
factor the augmentation into $\LG\to \LGamma\to \Z_p$ using (\ref{groupsurj})
and get: 
\begin{equation}\label{firstfact}
H^i(\Oh_K[1/S], \LG\otimes_{\Z_p}T(j))\to 
H^i(\Oh_K[1/S], \LGamma\otimes_{\Z_p}T(j)).
\end{equation}
The analysis of the cokernel of the
twisting map $Tw_{T(j-1)}$ will proceed in two
steps. The first is to investigate the cokernel of (\ref{firstfact}).
The second step treats then the cokernel of the map induced by the
augmentation $\LGamma\to \Z_p$:
\begin{equation}\label{secondfact}
H^i(\Oh_K[1/S], \LGamma\otimes_{\Z_p}T(j))\to 
H^i(\Oh_K[1/S], T(j)).
\end{equation}
\begin{lemma}\label{specseq} There is a spectral sequence
\begin{multline*}
E_2^{r,s}=\Tor_{r}^{\LG}\left( H^s(\Oh_K[1/S], \LG\otimes_{\Z_p}T(j)), \LGamma\right)\\
\Rightarrow H^{s-r}(\Oh_K[1/S],\LGamma\otimes_{\Z_p} T(j)).
\end{multline*}
\end{lemma}
\bew The projection formula (see e.g. \cite[Exercise 10.8.3]{We})
in the derived category gives
\[
R\Gamma(\Oh_K[1/S],\LG\otimes_{\Z_p}T(j)){\otimes}^{\La}_{\LG}\LGamma\isom
R\Gamma(\Oh_K[1/S],\LGamma\otimes_{\Z_p}T(j)).
\]
Taking cohomology gives the desired spectral sequence.
\bewende
\begin{cor}\label{surjobstr} There is an exact sequence
\begin{multline*}
H^1(\Oh_K[1/S], \LG(1))\otimes_{\Z_p} T(j-1)\xrightarrow{Tw_{T(j-1)}} H^1(\Oh_K[1/S], \LGamma\otimes_{\Z_p} T(j)) \\
\to \Tor_{1}^{\LG}\left( H^2(\Oh_K[1/S], \LG\otimes_{\Z_p}T(j)), \LGamma\right).
\end{multline*}
\end{cor}
\bew As $p>2$ the cohomological $p$-dimension of $\Oh_K[1/S]$ is $2$ 
 and
the result follows from the spectral sequence and the fact that the
twisting map factors through 
\[
H^1(\Oh_K[1/S], \LG\otimes_{\Z_p} T(j))\to 
H^1(\Oh_K[1/S], \LG\otimes_{\Z_p} T(j))\otimes_{\LG}\LGamma.
\]
\bewende
\begin{lemma}\label{torident} The canonical isomorphism $\LG{\otimes}_{\LH}\Z_p{\isom} \LGamma$
induces isomorphisms for all $r$:
\begin{multline*}
\Tor_{r}^{\LG}\left( H^2(\Oh_K[1/S], \LG\otimes_{\Z_p}T(j)), \LGamma\right)\isom \\
\Tor_{r}^{\LH}\left( H^2(\Oh_K[1/S], \LG\otimes_{\Z_p}T(j)), \Z_p\right).
\end{multline*}
In particular, one gets from corollary \ref{surjobstr} an exact sequence
\begin{multline*}
H^1(\Oh_K[1/S], \LG(1))\otimes_{\Z_p} T(j-1)\to  
H^1(\Oh_K[1/S], \LGamma\otimes_{\Z_p} T(j))\to \\
\to \Tor_{1}^{\LH}\left( H^2(\Oh_K[1/S], \LG\otimes_{\Z_p}T(j)), \Z_p\right).
\end{multline*}
\end{lemma}
\bew The isomorphism $\LG{\otimes}_{\LH}\Z_p{\isom} \LGamma$ can be checked at finite level as
$\Z_p$ is a finitely generated $\LH$-module. Then $\Z_p[G_n]\otimes_{\Z_p[H_n]}\Z_p\isom \Z_p[G_n/H_n]$ and
the claim is obvious. In particular, for finitely generated $\LG$-modules
$M$ is the functor $M\mapsto M\otimes_{\LG}\LGamma$ 
isomorphic to $M\mapsto M\otimes_{\LH}\Z_p$. 
\bewende
From this lemma and the factorization of the twisting map 
it is clear that the cokernel of the twisting map
is controlled by 
\[
\Tor_{1}^{\LH}\left( H^2(\Oh_K[1/S], \LG\otimes_{\Z_p}T(j)), \Z_p\right)
\] 
and
by $\Tor_{1}^{\LGamma}\left( H^2(\Oh_K[1/S], \LGamma\otimes_{\Z_p}T(j)), \Z_p\right)$.
To say something about these groups we need some results of Coates and
Sujatha on the finiteness of the $H^2$'s involved.
\subsection{Finiteness conditions for $H^2$}
This section contains only slight modifications of 
results of Coates and Sujatha \cite{C-S}. We thank them very much for making
these results available to us before their publication. One should
also compare this section with the appendix B in Perrin-Riou \cite{P-R} prop. B.2.

Let $\Lcyc$ (resp. $L_\infty$) be the maximal unramified abelian $p$-extension
of $\Kcyc$ (resp. $K_\infty$), 
in which every prime above $p$ splits completely.
\begin{prop}\label{CS}(Coates-Sujatha \cite{C-S}) Assume that $\Gh=\Gal(K_\infty/K)$ 
is a pro-$p$-group, 
then the following conditions are equivalent:
\begin{itemize}
\item[i)] $\Gal(\Lcyc / \Kcyc)$ is a finitely generated $\Z_p$-module
\item[ii)] $\Gal(L_\infty / K_\infty)$ is a finitely generated $\LH$-module
\item[iii)] $H^2(\Oh_K[1/S], \LGamma\otimes_{\Z_p}T)$ is 
a finitely generated $\Z_p$-module
\item[iv)] $H^2(\Oh_K[1/S], \LG\otimes_{\Z_p}T)$
is a finitely generated $\LH$-module.
\end{itemize}
In particular, if these equivalent
conditions are satisfied, the $\LGamma$-module 
$H^2(\Oh_K[1/S], \LGamma\otimes_{\Z_p}T)$ is torsion, 
i.e., the weak Leopold conjecture is true.
\end{prop}
{\bf Remarks:} The first condition is equivalent to the
famous $\mu=0$ conjecture of Iwasawa (see \cite{NSW} Ch. XI thm. 11.3.18), 
which is known to be true for $K/ \Q$ abelian.

Note also that the statements $i)$ and $ii)$ in the proposition
are independent of the Galois representation $T$.

\bew The proof of the proposition can be found in
Coates and Sujatha \cite{C-S} in the case 
of the Tate module for an elliptic curve. The case 
for an arbitrary Galois representation $T$ is the same.
More precisely, the equivalence $i)\iff ii)$ is lemma
3.7., $i)\iff iii)$ is thm. 3.4. in loc. cit. To prove $iii)\iff iv)$,
we have from the spectral sequence in lemma \ref{specseq}
and the vanishing of  {\'e}tale cohomology for $s>2$ that
\[
H^2(\Oh_K[1/S], \LG\otimes_{\Z_p}T)\otimes_{\LH}\Z_p\isom
H^2(\Oh_K[1/S], \LGamma\otimes_{\Z_p}T)
\]
and the claim follows from Nakayama's lemma.
\bewende

\noindent{\bf Example:} Let  $E/\Q$ be an elliptic curve over $\Q$ and
$F/\Q$ be an abelian extension such that $E_{p^\infty}(F)\neq 0$. Then
it is easy to see (cf. \cite{C-S} cor. 3.6.) that $F(E{_p^\infty})/F(\mu_p)$
is a pro-$p$ extension. Thus these elliptic curves provide examples
where the above proposition \ref{CS} applies. More specific 
examples are $E: y^2+xy=x^3-x-1$ and $F=\Q(\mu_7)$ or
$E:  y^2+xy=x^3-3x-3$ and $F=\Q(\mu_5)$ (see loc. cit. 4.7. and 4.8.).

\subsection{{\'E}tale cohomology classes as twists of units}\label{keysection}
Recall that $\Lcyc$ is the maximal unramified abelian $p$-extension
of $\Kcyc$, in which every prime above $p$ splits completely,
and that we have an isomorphism $H^1(\Oh_K[1/S], \LG(1)) \cong
\prolim_n\Oh_{K_n}[1/S]^*\otimes\Z_p$ 
by Proposition \ref{twistisom} and (\ref{Kummer}).
\begin{thm}\label{surjcond} Suppose  that $\Gal(\Lcyc / \Kcyc) $ is a
finitely generated $\Z_p$-module, then there exists a $J\in \N$ such that for 
all $j\geq J$ the twisting map 
\[
H^1(\Oh_K[1/S], \LG(1))\otimes_{\Z_p}T(j-1)\xrightarrow{\Tw_{T(j-1)}} 
H^1(\Oh_K[1/S], T(j))
\]
has finite cokernel. In particular,
for $j\geq J$ all elements in $H^1(\Oh_K[1/S], V(j))$ (where $V=T\otimes_{\Z_p}\Q_p$ as before) 
are ``twists'' of norm compatible units in 
\[
H^1(\Oh_K[1/S], \LG(1))=\prolim_n H^1(\Oh_{K_n}[1/S], \Z_p(1))
\]
with a basis in the lattice $T(j-1)$.
\end{thm}
{\bf Remark:} The choice of the twist $1$ in $H^1(\Oh_K[1/S], \LG(1))$
and hence of the 
group of norm compatible units instead of any other twist 
is just for esthetic reasons.
For the application to Euler systems and the construction of $p$-adic
L-functions the units are certainly the most interesting case. In particular,
we see this theorem as a strong confirmation of the philosophy explained
in \cite{Hu-Ki2}, that all $p$-adic properties of motives in connection
with L-values should be encoded in the associated tower of number fields.

It is an interesting question to investigate $H^2(\Oh_K[1/S], T(j))$
with the above methods and to compare this with the results by McCallum
and Sharifi \cite{SMcC}.

\bew  It follows from proposition \ref{CS} that under our conditions
\[
\Tor_1^{\LH}\left(H^2(\Oh_K[1/S], \LG\otimes_{\Z_p}T(j)),\Z_p\right)
\]
is a finitely generated $\Z_p$-module. Indeed $H^2(\Oh_K[1/S], \LG\otimes_{\Z_p}T(j))$ is
a finitely generated $\LH$-module and thus the groups 
\[
\Tor_r^{\LH}\left(H^2(\Oh_K[1/S], \LG\otimes_{\Z_p}T(j)),\Z_p\right)
\]
are also finitely generated $\Z_p$-modules by standard homological algebra.
The exact sequence \ref{surjobstr} implies that the
cokernel of 
\[
H^1(\Oh_K[1/S], \LG\otimes_{\Z_p}T(j))\otimes_{\LG}\LGamma\to H^1(\Oh_K[1/S], \LGamma \otimes_{\Z_p}T(j))
\]
is a $\LGamma$-module, say $M(j)$, which is finitely generated as a
$\Z_p$-module, hence torsion as $\LGamma$-module. 
By the classification of torsion $\LGamma$-modules, the coinvariants 
of $M(j)\otimes_{\LGamma}\Z_p$ are finite for sufficiently big $j$. 
We get an exact sequence
\begin{multline*}
H^1(\Oh_K[1/S], \LG\otimes_{\Z_p}T(j))\otimes_{\LG}\Z_p\to\\
H^1(\Oh_K[1/S], \LGamma \otimes_{\Z_p}T(j))\otimes_{\LGamma}\Z_p
\to M(j)\otimes_{\LGamma}\Z_p\to 0.
\end{multline*}
To get the twisting map we have to compose with the first map in the following
exact sequence (which is similar to Corollary \ref{surjobstr})
\begin{multline*}
H^1(\Oh_K[1/S], \LGamma \otimes_{\Z_p}T(j))\otimes_{\LGamma}\Z_p\to\\ 
H^1(\Oh_K[1/S], T(j))\to 
\Tor_1^{\LGamma}\left(H^2(\Oh_K[1/S], \LGamma\otimes_{\Z_p}T(j)),\Z_p\right).
\end{multline*}
By our condition and proposition \ref{CS} $H^2(\Oh_K[1/S], \LGamma\otimes_{\Z_p}T(j))$ is also a finitely 
generated $\Z_p$-module and thus $\LGamma$-torsion. 
As $\Gamma$ is cyclic ($\Gh$ and hence $\Gamma$ is pro-$p$),
the $\Tor_1^{\LGamma}$ term identifies with
the $\Gamma$-invariants of $H^2(\Oh_K[1/S], \LGamma\otimes_{\Z_p}T(j))$.
 Again for $j$ big enough these are finite.
\bewende
{\bf Remark:} In fact, if $M$ is the $p$-adic realization 
of a motive (say $M=H^n(X\times_K\Kbar, \Q_p)$ for $X/K$ smooth, projective), 
one should expect that the $\LGamma$-module
\[
\Tor_1^{\LG}\left(H^2(\Oh_K[1/S], \LG\otimes_{\Z_p}T(j)),\LGamma\right)\otimes_{\LGamma}\Z_p
\]
is finite for all $j\geq n+1$. Compare this with Jannsen's
conjecture \ref{Jannsen} about the vanishing of 
$H^2(\Oh_K[1/S],H^n(X\times_K\Kbar, \Q_p(j)))$ for $j\geq n+1$.

\section{Twisting for motivic cohomology with $p$-adic coefficients}
In this section  $X$ will always be a 
smooth and projective scheme over $D=\Oh_K[1/S]$.

The goal in this section is to study the twisting map in
motivic cohomology with finite coefficients. The general
assumption is that the ``Bloch-Kato-conjecture'' for motivic
cohomology holds 
as announced by Voevodsky in \cite{Vo2}
(do not confuse this with the Tamagawa number conjecture). This implies, using 
the Beilinson-Lichtenbaum conjecture, that we have to deal with
{\'e}tale cohomology of $X$.

\subsection{Review of motivic cohomology with finite coefficients
over Dedekind domains}
For a variety $X$ smooth over a Dedekind ring $D$, we define
motivic cohomology groups as the hypercohomology of Bloch's cycle complex
$\Z(j)$. As usual, $\Delta^s_D:=\Spec(D[t_0,...,t_s]/\sum_i t_i -1)$ 
denotes the standard algebraic $s$-simplex.

For a variety $X$ smooth over a Dedekind ring $D$,
let $z^j(X,i)$ be the free abelian group on closed integral subschemes
of codimension $j$ on $X \times_{D} \Delta^i_D$ which intersect all faces properly.
The associated complex of presheaves (with $z^j(X,2j-i)$ in degree $i$)
is denoted $\Z(j)$, and
$\Z/n(j):=\Z(j)\otimes^{\mathbb L} \Z/n$. The complex $\Z(j)$ (and thus also 
$\Z/n(j)$) is a complex of sheaves for the {\'e}tale topology \cite[Lemma 3.1]{Ge},
and we write $\Z/n(j)_{et}$ resp. $\Z/n(j)_{Zar}$ when considering it as a complex 
of {\'e}tale resp. Zariski sheaves.
\begin{defn}\label{defmotcohom}(compare \cite[p. 5]{Ge})
The {\em motivic cohomology}
of $X$ is the hypercohomology  
\begin{equation}
H^i_{mot}(X,\Z/n(j)):=H^i(X,\Z/n(j)_{Zar}).
\end{equation}
\end{defn}

Calling this motivic cohomology is justified by Voevodsky's \cite{Vo1} theorem
$Hom_{DM^{eff,-}(K)}(M(X),\Z(j)[i])=:H_{mot}^i(X,\Z(j))\cong
CH^j(X,2j-i)$ if $D=K$ is a field. 
In this case, higher Chow groups are defined by taking just 
cohomology and not hypercohomology. By \cite[Theorem 3.2]{Ge} 
both definitions coincide not only over a field but still
if the base $D$ is a discrete valuation ring.

Observe \cite[section 3]{Ge} that $H^i_{mot}$ is covariant for 
proper maps (with degree shift) and contravariant for flat maps. The latter 
applies in particular to the structural morphisms
$p_n:X_n \to D_n$.

The {\'e}tale cycle class  $cl$ factors
through the {\'e}tale sheafification $\Z/n(j)_{et}$ via the
map $\Z/n(j)_{Zar} \to {\bf R}\pi_*\Z/n(j)_{et}$ 
induced by the morphism of sites $\pi:(Sm/D)_{et} \to (Sm/D)_{Zar}$.

For us the most important consequence of the Bloch-Kato conjecture
is the truth of the Beilinson-Lichtenbaum conjecture:

\begin{thm}(Geisser \cite[Theorem 1.2 (2)(4)]{Ge})\label{Geisser}
Assume that X is a smooth scheme over a Dedekind domain
$D$ with $n \in D^{\times}$ and that the ``Bloch-Kato-conjecture''
holds.
\begin{itemize}
\item[1)] For all $i$ and $j$ there is an isomorphism
\[
H^i(X,\Z/n(j)_{et})\isom H^i(X, \Z/n(j))
\]
of the {\'e}tale hypercohomology of $\Z/n(j)_{et}$ with the {\'e}tale
cohomology.
\item[2)]
The {\'e}tale cycle class map
induces isomorphisms for $0\leq i\leq j$
\[ H^{i}_{mot}(X,\Z/n(j)) \cong H^{i}(X,\Z/n(j)) \]
of motivic with {\'e}tale cohomology.
\end{itemize}
\end{thm}

%\begin{defn} The $p$-adic motivic cohomology is
%\[
%H^i_{mot}(X,\Z_p(j)):= \prolim_n H^i_{mot}(X,\Z/p^n(j)).
%\]
%\end{defn}

\subsection{The geometric twisting map}
We are going to define a geometric twisting map, which will allow 
to relate our results for {\'e}tale cohomology with motivic cohomology.
The main difficulty is that the cup-product is not compatible with
corestriction maps. We use the compatibility of the Hochschild-Serre
spectral sequence with cup-product to overcome this and to reduce 
to an observation due to Soul{\'e}.

In this section we consider $X\to \Spec\Oh_K[1/S]$ smooth and proper.
We denote by $\Xbar:=X\times \Kbar$ and
let $T:= H^i(\Xbar, \Z_p)$. This is a Galois-module and
finitely generated $\Z_p$-module.
As in section \ref{sectiontwist} this defines a 
tower of number fields $K_n$ and a $p$-adic
Lie group $\Gh:=\Gal(K_\infty /K)$. Let   
\[
X_n:=X\times_{\Oh_K[1/S]}\Oh_{K_n}[1/S].
\]
To construct a twisting map for motivic cohomology we use the pull-back
\[
H^1(\Oh_{K_n}[1/S],\Z/p^n(1))\to H^1(X_n,\Z/p^n(1))
\]
and the cup-product with $H^{i}(X_n,\Z/p^n(j-1))$. This produces 
elements in $H^{i+1}(X_n,\Z/p^n(j))$. The Hochschild-Serre spectral sequence
\[
E_2^{p,q}=H^p(\Oh_{K_n}[1/S],H^q(\Xbar, \Z/p^n(j)))\Rightarrow
H^{p+q}(X_n,\Z/p^n(j))
\]
allows to relate these elements with 
$H^1(\Oh_{K_n}[1/S],H^i(\Xbar, \Z/p^n(j)))$ as follows:
Let 
\begin{multline*}
H^{i+1}(X_n,\Z/p^n(j))^0:=\\
\ker\left( H^{i+1}(X_n,\Z/p^n(j))\xrightarrow{\gamma}
H^0(\Oh_{K_n}[1/S],H^{i+1}(\Xbar, \Z/p^n(j)))\right)
\end{multline*}
be the kernel of the edge morphism $\gamma$. As the Hochschild-Serre 
spectral sequence is compatible with cup-products, we get a map
\begin{equation}\label{cupprod}
H^1(\Oh_{K_n}[1/S],\Z/p^n(1))\times H^{i}(X_n,\Z/p^n(j-1))\to 
H^{i+1}(X_n,\Z/p^n(j))^0.
\end{equation}
As $E_2^{1,i}=E_\infty^{1,i}$ we have a surjection,
\begin{equation}\label{suredge}
H^{i+1}(X_n,\Z/p^n(j))^0\to
H^1(\Oh_{K_n}[1/S],H^i(\Xbar, \Z/p^n(j))),
\end{equation}
which we compose
with the corestriction map
\begin{equation}\label{cores}
H^1(\Oh_{K_n}[1/S],H^i(\Xbar, \Z/p^n(j)))\xrightarrow{\cores}
H^1(\Oh_{K}[1/S],H^i(\Xbar, \Z/p^n(j))).
\end{equation}
If we compose the cup-product in (\ref{cupprod})
with this composition and using again the compatibility of the 
spectral sequence with products
we see that the cup-product has to factor through
the edge morphism 
\begin{equation}\label{edge}
H^{i}(X_n,\Z/p^n(j-1))\xrightarrow{\gamma}
H^0(\Oh_{K_n}[1/S],H^{i}(\Xbar, \Z/p^n(j-1))).
\end{equation}
Thus we get
\begin{multline}\label{geometrictwist}
H^1(\Oh_{K_n}[1/S],\Z/p^n(1))\times 
H^0(\Oh_{K_n}[1/S],H^{i}(\Xbar,\Z/p^n(j-1)))\to\\ 
\to H^1(\Oh_{K}[1/S],H^i(\Xbar, \Z/p^n(j))).
\end{multline}
It is an important observation by Soul{\'e} that, although
the cup-product in general  is not compatible with corestriction,
the map in (\ref{geometrictwist}) is compatible.
To formulate the observation of Soul{\'e} properly we need:
\begin{defn}
A sequence 
of elements $\alpha_n \in H^r(X_n,\Z/p^n)$ is {\em norm compatible}
if the restriction of the coefficients modulo $\Z/p^{n-1}$
of $cores(\alpha_n)$ is $\alpha_{n-1}$ for all $n \geq 2$.
The sequence $\{\alpha_n\}$ is {\em reduction compatible} if 
the reduction modulo $p^{n-1}$ of $\alpha_{n}$ 
is the pull-back of $\alpha_{n-1}$ for all $n \geq 2$. 
\end{defn}
Note that the elements $\{\alpha_n\}$ of any $\Z_p$-lattice 
$T\subset H^i(\Xbar, \Q_p(j)))$ are reduction compatible.
Soul{\'e} proves the following: 
\begin{lemma}\label{NRN}
If the sequence $\{\alpha_n\}$ is norm compatible and $\{\beta_n\}$ 
is reduction compatible, then $\{\alpha_n \cup \beta_n\}$ is 
norm compatible.
\end{lemma}
\proof
This is just the projection formula, see \cite[Lemma 1.4]{So2}.
\qed
Taking the projective limit over $n$ in (\ref{geometrictwist}) gives:
\begin{cor} Cup-product gives a map 
\begin{multline*}
H^1(\Oh_{K}[1/S],\LG(1))\otimes 
H^{i}(\Xbar,\Z_p(j-1))\to H^1(\Oh_{K}[1/S],H^i(\Xbar, \Z_p(j))).
\end{multline*}
\end{cor}
To apply theorem \ref{surjcond} in this situation, one has to 
choose a lattice $T\subset H^i(\Xbar, \Z_p)$, which is stable under the
action of the Galois group. Observing that for a finite Galois module
$M$ the cohomology $ H^1(\Oh_{K}[1/S],M)$ is finite, we get:
\begin{cor}\label{twistcor}
Under the condition of theorem \ref{surjcond} there is 
an integer $m$ such that for all $n\geq m$ the cokernel of
the cup-product map in (\ref{cupprod}) composed with the maps
in (\ref{suredge}) and (\ref{cores})
\begin{multline*}
H^1_{mot}(\Oh_{K_n}[1/S],\Z/p^n(1))\times H^{i}_{mot}(X_n,\Z/p^n(j-1))\to\\
\to  H^1(\Oh_{K}[1/S],H^i(\Xbar, \Z/p^n(j)))
\end{multline*}
is annihilated by $p^m$ for all $j\geq J$.
\end{cor}
{\bf Remarks:} a) More generally,
the above construction is possible for  any theory $A^*$, which is 
covariant for proper maps, contravariant
for flat maps and satisfying the projection formula 
$f_*(a \cup f^*(b))=f_*(a) \cup b$ for flat proper 
(or at least finite {\'e}tale) maps $f$. We explain the case of
K-theory with finite coefficients in appendix \ref{Kappendix}.

b) Soul{\'e} applies the above  construction
to get non-torsion elements
in the $K$-groups of rings of integers or elliptic curves with
complex multiplication. In these cases the schemes $X_n$ are base changes
of $X$ to the ring $\Oh_{K_n}[1/S]$, where $K_n$ is defined by adjoining
$p^n$-th roots of unity or $p^n$-th division points of the elliptic curve.
The towers of fields are in these cases abelian. It is shown in the cyclotomic
case in \cite{Hu-Wi} and \cite{Hu-Ki1} (with another method)
and in the case of CM-elliptic curves in \cite{Ki} that these twisted
elements are in fact motivic, i.e., are in the image of motivic cohomology.

\subsection{Compatibility of cup products in motivic and {\'e}tale cohomology}
\label{compat}
The aim of this technical section is to establish the compatibility of
cup products in {\'e}tale and motivic cohomology. The problem is
that the cup product for motivic cohomology over Dedekind rings
is only defined if one factor consists of equi-dimensional cycles 
(see definition \ref{defprod} below).
We will show that we have 
a commutative diagram 
\begin{equation}\label{diagram}
\xymatrix{
(D_n^{\times} \otimes \Z/p^n) \times H_{mot}^i(X_n,\Z/p^n(j-1))
\ar[r]^(.61){\cup_{mot}  \circ {\phi}\times \id} 
\ar[d]^{{\phi}\times cl} 
\ar@{}[dr]^{} & H_{mot}^{i+1}(X_n,\Z/p^n(j)) 
\ar[d]^{cl}  \\
H^1(D_n,\Z/p^n(1)) \times H^i(X_n,\Z/p^n(j-1))
\ar[r]^(.6){\cup} \ar[d]^{1 \times  \gamma} & H^{i+1}(X_n,\Z/p^n(j)) \ar[d]^{\tilde{\gamma}}\\
H^1(D_n,\Z/p^n(1)) \times H^i((\Xbar, \Z/p^n(j-1)) 
\ar[r]^(.6){Tw} & H^1(D_n,H^i(\Xbar, \Z/p^n(j)) 
}
\end{equation}
where $D_n:=Spec(\Oh_{K_n}[1/S])$, 
$\gamma$ and $\tilde{\gamma}$ are the edge maps as before,
the vertical arrows $cl$ are {\'e}tale cycle class maps
and the cup product $\cup_{mot}$ as well as the map ${\phi}$
are defined below.
The commutativity of the lower square follows from
lemma \ref{twistiscup} and the compatiblity of the Hochschild-Serre
spectral sequence with cup products.
The commutativity of the upper square of (\ref{diagram}) is discussed below;
this generalizes the classical result for the
usual cycle class map (see e. g. \cite[Proposition VI.9.5]{Mi}).
Recall that by \cite[Corollary 4.3]{Su} we have an isomorphism
$H_{mot}^i(\Xbar,\Z/p^n(j)) \cong H^i(\Xbar,\Z/p^n(j))$.

Recall from \cite[section 1.7]{Le} that an irreducible scheme
$Z\to D$ is {\em equi-dimensional} if it is dominant over $D$.
The relative dimension $\dim_{D}Z$ 
is then defined to be the dimension of the
generic fibre.
Now we define the relative higher Chow group complex
for our smooth $X\to D$ as follows: $z^{j}(X/D,p)$ to be
the free abelian group generated by irreducible closed
subsets $Z\subset X \times_D \Delta^p_D$, such that for each face $F$ of 
$\Delta^p_D$ the irreducible components $Z'$ of $Z\cap (X\times F)$
are equi-dimensional over $D$ and $\dim_{D}Z'=\dim_{D}F+d-j$.
Note that  we have an inclusion of complexes $z^j(X/D,*) \subset z^j(X,*)$.
We define equi-dimensional motivic cohomology 
$H^i_{mot}(X/D, \Z(j))$ to be the Zariski
hypercohomology  of the complex which has in degree $i$
the Zariski sheafification of the presheaf $U\mapsto z^j(U/D,2j-i)$. 
To define $H^i_{mot}(X/D, \Z/p^n(j))$ we use the same complex tensored with 
$\otimes^{\mathbb L}\Z/p^n$.

The units $D_n^{\times}$ we use for twisting are all
equi-dimensional:

\begin{lemma} 
The map 
\[
\phi:D_n^{\times}{\to} H^1_{mot}(D_n,\Z(1))
\]
induced by sending $u\neq 1\in D_n^{\times}$ to the graph
of the rational map 
\[
\left(\frac{1}{1-u},\frac{u}{u-1}\right):\Spec D_n\to \Delta^1_D
\]
(i.e. to a cycle in $D_n\times_D\Delta^1_D$)
factors through $H^1_{mot}(D_n/D,\Z(1))$. The induced map
\[
D_n^{\times}\otimes\Z/p^n\to H^1_{mot}(D_n/D,\Z/p^n(1))
\]
is injective.
\end{lemma}
\bew 
In \cite[Lemma 11.2]{Le} Levine constructs a map
$D_n^{\times}{\to} CH^1(D_n,1)$ using the graph of 
$\left(\frac{1}{1-u},\frac{u}{u-1}\right)$. Together with the
natural map $CH^1(D_n,1)\to H^1_{mot}(D_n,\Z(1))$ this defines $\phi$
and hence a map 
\[
D_n^{\times}\otimes\Z/p^n\to H^1_{mot}(D_n,\Z/p^n(1)).
\]
If we compose this with the isomorphism
in \ref{Geisser}, we get a map
\[
D_n^{\times}\otimes\Z/p^n\to H^1_{et}(D_n,\Z/p^n(1)),
\]
which is obviously (reduce to the case of a field) the map 
induced by the Kummer sequence, hence injective.
It remains to show that the map factors through $H^1_{mot}(D_n/D,\Z/p^n(1))$.
As the graph of $\left(\frac{1}{1-u},\frac{u}{u-1}\right)$ is an
equi-dimensional cycle, this follows from the definition.
\bewende
Now we define the upper horizontal map $\cup_{mot}$ of (\ref{diagram}).
\begin{defn}\label{defprod}
For $f: X_n \to Spec(D_n)$ smooth,
we define 
\[ \cup_{mot}:
H_{mot}^1(D_n/D,\Z/p^n(1)) \times H_{mot}^i(X_n,\Z/p^n(j-1))
\to H_{mot}^{i+1}(X_n,\Z/p^n(j)) \]
as the composition
\begin{multline*}
H_{mot}^1(D_n/D,\Z/p^n(1)) \times H_{mot}^i(X_n,\Z/p^n(j-1))
\xrightarrow{(f\times \id)^*\circ\cup^{r,s}_{D_n/D,X_n}}\\ 
\to H_{mot}^{i+1}(X_n \times_D D_n,\Z/p^n(j))
\xrightarrow{(id,p_n)^*}H_{mot}^{i+1}(X_n,\Z/p^n(j)).
\end{multline*}
\end{defn}
Here $\cup^{r,s}_{D_n/D,X}:z^s(D_n/D,*) \otimes z^r(X) \to z^{r+s}(X\times_{D} D_n)$ 
is the exterior product with integral coefficients defined by
Levine \cite[section 8]{Le}.
The product of the complexes of presheaves 
induces a product of complexes of sheaves and
(using Godement resolutions as in \cite{GL}) 
on the hypercohomology groups.

Now we return to the commutativity of (\ref{diagram}). By definition
of the twisting map at finite level in \ref{twistvariant}, it is 
enough to consider the
diagram
\[
\xymatrix{
H^1_{mot}(D_n/D,\Z/p^n(1))\times H_{mot}^i(X_n,\Z/p^n(j-1))
\ar[r]^(.6){\cup_{mot}} \ar[d]^{cl} & 
H_{mot}^{i+1}(X_n \times_D D_n,\Z/p^n(j)) 
\ar[d]^{cl}  \\
H^1(D_n,\Z/p^n(1)) \times H^i(X_n,\Z/p^n(j-1))
\ar[r]^(.6){\cup}  & H^{i+1}(X_n \times_D D_n,\Z/p^n(j)). \\
}
\]
As pointed out in \cite[p. 13]{Ge}, the proof of
\cite[Proposition 4.7]{GL} for the commutativity of the corresponding 
diagram of varieties over fields
carries over to Dedekind domains.
The argument in the proof of \cite[Proposition 4.7]{GL}
that $\cup$ equals the product $\cup'$ of {\it loc. cit.} constructed
in a way compatible with $\cup_{D_n/D,X}$ is still valid
over Dedekind domains. Hence the commutativity of (\ref{diagram}).

%We also have the following projection formula
%for motivic cohomology over Dedekind rings,
%which shows that lemma \ref{NRN} applies
%to motivic cohomology.

%\begin{lemma}\label{projmot}
%Let $a \in H_{mot}^1(D_n,\Z/p^n(j-1))$,
%$b \in H_{mot}^i(X_{n-1},\Z/p^n(j))$
%and $f_n:X_n \to X_{n-1}$ be the finite {\'e}tale map 
%induced by $g_n:D_n \to D_{n-1}$.
%Then we have 
%$$(f_n)_*(a \cup_{mot} f_n^*b)=(g_n)_*a \cup_{mot} b$$.
%\end{lemma}
%\bew
%This follows by the definition of $\cup_{mot}$ and the base 
%change fromula \cite[Proposition 2.3.4]{SV} which implies that 
%$(f_{n-1})_* \circ (p_n,f_n)^* = (p_{n-1},id)^* \circ (q_n 
%\times id)_*$, and similarly after taking the fibre product over $D$ 
%with $\Delta^1_D \times_D \Delta^{2j-2-i}_D$.
%\bewende 

%\rem It is not known if the groups of Definition \ref{defmotcohom}
%are isomorphic to the groups represented
%by motivic Eilenberg-Mac Lane spaces in the motivic
%homotopy category over the base $S=Spec(\Oh)$ of Morel and 
%Voevodsky \cite{MV}.
% Die Referenzen sind
% \bibitem[Bl]{Bl} S. Bloch: Algebraic Cycles and Higher $K$-theory,
% Adv. in Math. 61 (1986), 267-304
% \bibitem[Fu]{Fu} K. Fujiwara: A proof of the absolute purity 
% conjecture (after Gabber). Algebraic geometry 2000, Azumino (Hotaka), 153--183, 
% Adv. Stud. Pure Math., 36, Math. Soc. Japan, Tokyo, 2002.

\begin{appendix}

\section{Twisting in $p$-adic $K$-theory}\label{Kappendix} 
In this appendix, we will reinterpret our results
in terms of $p$-adic $K$-theory.

As usual, we can define $K$-theory with coefficients
of the exact category $Vect(X)$ of vector bundles 
on $X$ using Quillen's Q-construction and homotopy groups with 
finite coefficients: 
\begin{defn}
Let $$ K_m(X,\Z/q):=\pi_m(\Omega BQVect(X),\Z/q)$$
for $r>0$ and $K_0(X,\Z/q):=K_0(X)/q$.
Moreover, we set
$$ K_r(X,\Z_p):=\prolim_n \, K_r(X,\Z/p^n).$$ and define
$K_r(X,\Q_p):=K_r(X,\Z_p)\otimes_{Z_p}\Q_p$.
\end{defn}

Here we use that we have maps $\eps_n:K_r(X,\Z/p^n) \to
K_r(X,\Z/p^{n-1})$ given by reduction of coefficients.
Applying $\prolim$ to the short exact sequence 
$$ 0 \to K_r(X)/p^n \to K_r(X,\Z/p^n) \to _{p^n}K_{r-1}(X) \to 0$$
shows that $rk_{\Z}K_r(X)=rk_{\Z_p}K_r(X,\Z_p)$, 
provided the groups $K_r(X)$ and $K_{r-1}(X)$ are 
finitely generated as generally conjectured 
(``Bass conjecture'') and proved
if $X=\Spec(\Oh_K)$ by Quillen \cite{Qu}.

\medskip

%Die folgende Vermutung fliegt wohl raus.
%We will see below that for $j > > 0$ the map
%$\gamma$ induces a map $\tilde{\gamma}:H^n(X_r,\Z/p^r(j-1)) \to
%K_{2j-n-2}(X_r,\Z/p^r)^{(j-1)}$.
%\begin{conj}\label{main}
%Let $d=dim(X \times_{{\Oh}_K} K)$ (or equivalently
%$d=dim(X)-1$), and assume that $m\geq (8/3)(d+2)(d+3)
%(d+4)$ and $p$ an odd prime . Then
%$$rk_{{\Z}_p}K_m(X,{\Z}_p)^{(i)} \geq n_{i,2i-m}$$
%and the $n_{i,2i-m}$ generators of $K_m(X,{\Z}_p)^{(i)}
%\otimes \Q_p$
%can be expressed as $r_p(N(\beta \cup \tilde{\gamma}(\alpha)))$
%where $(\alpha)$ is induced by some $a \in H^n(X\times_{{\Oh}_K[1/s]}\bar{K},
%\Z_p(j))$ and $(\beta)$ is a family of elements
%in $K_1({\Oh}_K[1/S]_r,\Z/p^r)^{(1)}$ having property $N$.  
%\end{conj}

We assume as before that $X$ is smooth over $\Oh_K[1/S]$,
of relative dimension $d$. 
Adams operations carry over to finite coefficients
and their eigenspaces are denoted by $K(X,\Z/p^n)^{(j)}$ as usual.
By \cite[Proposition 6]{So1.5} the transfer maps $(f_n)_*$ respect
these eigenspace decomposition (the hypothesis of {\it loc. cit.}
is satisfied as the field extension $K_{n}/K_{n-1}$
is finite).

%The proof of Conjecture  \ref{main} will consist in comparing this pairing with the 
%corresponding one in {\'e}tale cohomology and then using the results of the previous 
%sections.

\medskip

Thomason constructs an algebraic Bott
element $\beta \in K_2(X,\Z/p^n)$ and proves that
there is an isomorphism $K_*(X,\Z/p^n)[\beta^{-1}] 
\stackrel{\cong}{\to} K^{et}_*(X,\Z/p^n)$
\cite[Theorem 4.11]{Th1}, that 
$\phi_j:K_j(X,\Z/p^n) \to K^{et}_j(X,\Z/p^n)$
is an epimorphism if $j \geq N$ and $\beta ^N$ annihilates $ker(\phi_j)$
for all $j \geq 0$, where $N=2/3(d+2)(d+3)(d+4)$ \cite[Corollary 3.6]{Th2}.
Multyplying the short exact (for $j \geq 2N$) sequence
$\ker(\phi_j) \to K_j(X,\Z/p^n) \to K^{et}_j(X,\Z/p^n)$ 
with $\beta^N$ and applying the snake lemma,
we get a splitting $K_{j+2N}(X,\Z/p^n) \to ker(\phi_{j+2N})$
and thus {\'e}tale $K$-theory is a natural direct 
summand of $K$-theory in these degrees.
So if $2j-i-2 \geq (8/3)(d+2)(d+3)(d+4)$, we obtain a pairing
$$K^{et}_1({\Oh}_{K_n}[1/S],\Z/p^n) \times K^{et}_{2j-i-2}(X_n,\Z/p^n) 
\to K^{et}_{2j-i-1}(X_n,\Z/p^n)$$
which is a direct summand of the corresponding pairing 
in algebraic $K$-theory with finite coefficients.
Concerning the first factor, we even have an isomorphism
between $K_1$ and $K_1^{et}$ by \cite[Proposition 8.2]{DF}.

\rem For $p=2$, the bounds for $j$ 
such that $K_j(X,\Z/p^n) \stackrel{\cong}{\to} K^{et}_j(X,\Z/p^n)$ have been 
improved by Kahn \cite[Theoprem 2]{Ka} provided $X$ is ``non-exceptional''. 
He shows that it is an isomorphism if $j \geq cd_2X - 1$
As he points out \cite[p. 104]{Ka},
these improved bounds will carry over to odd $p$ 
(without the non-exceptional restriction) assuming the Bloch-Kato conjecture 
for $K$ holds. 

The next step is to observe that the {\'e}tale
Atiyah-Hirzebruch spectral sequence degenerates
$(E_2=E_{\infty})$ provided 
$p>(cd_pX/2) +1$ where $cd_pX$ is the $p$-cohomological
dimension of $X$, which is at most $2d+3$
(see \cite[Expos{\'e} X]{SGA4}). 
Moreover, the Adams filtration on $K$-theory and 
the weight filtration on {\'e}tale cohomology coincide in a 
certain range \cite[Theorem 2]{So1}, so that
the left hand side of the above pairing for $K^{et}$ is 
isomorphic to $$H^1({\Oh}_{K_n}[1/S],\Z/p^n(1)) \otimes 
H^i(\Xbar,\Z/p^n(j-1))$$
provided $p \geq (j+cd_pX + 3)/2$.
As $H^i(\Xbar,\Z/p^n(j-1))$ is a trivial 
${\Oh}_{K_n}[1/S]$-sheaf, the twist of Definition \ref{twistingmap}
yields an isomorphism
$$H^1({\Oh}_{K_n}[1/S],\Z/p^n(1)) \otimes 
H^i(\Xbar ,\Z/p^m(j-1)) \cong
H^1({\Oh}_{K_n}[1/S],H^i(\Xbar,\Z/p^n(j))).$$
It is now possible to construct elements having property R and N for algebraic
$K$-theory, and to proceed as in the previous section.
The $p$-adic cycle class map has to be replaced by the 
$p$-adic regulator (take the inverse limit of 
\cite[Definition 2.22, Example 1.4.(iii)]{Gi})
\[
r_p:K_{2j-i-1}(X,\Z_p)^{(j)} \to H^{i+1}(X,\Q_p(j))).
\]

\section{Calculation of the inverse limit of Galois cohomology}\label{twistbew}
Here we give the proof of proposition \ref{twistisom}. Let $G:=G_S$ and
$H'\subset H\subset G$ subgroups defining $K_n$ and $K_m$,
so that $G/H'\isom G_m$ and $G/H\isom G_n$ (hence H/H' is finite).

By Shapiro's lemma we have
\[
H^i(\Oh_{K_n}[1/S], T)\isom H^i(\Oh_{K}[1/S],\Hom_H(G, T)),
\]
where $\Hom_H(G, T)$ denotes the continuous maps $f:G\to T$ such
that $f(hg)=hf(g)$. The group $G$ acts on this via $(gf)(x):=f(xg)$.
The corestriction on the left hand side 
\[
\cores:H^i(\Oh_{K_m}[1/S], T)\to H^i(\Oh_{K_n}[1/S], T)
\]
is induced on the right hand side by the map
\begin{align*}
\tr:\Hom_{H'}(G, T)&\to \Hom_H(G, T)\\
f&\mapsto\{g\mapsto \sum_{h\in H/H'}hf(h^{-1}g)\}
\end{align*}
A straightforward calculation shows that this is well-defined.
Consider now $T/p^n$ so that the $G$ action factors through $G_n$.
Define $^gf(x):= gf(g^{-1}x)$. Then because $H\subset G$ is
a normal subgroup this defines another $G$ action on $\Hom_H(G, T/p^n)$
where $H$ acts trivially.
Consider the $G$-isomorphism
\begin{align*}
\phi:\Hom_H(G, T/p^n)&\xrightarrow{\isom}\Z_p[G_n]\otimes_{\Z_p}T/p^n\\
f&\mapsto \sum_{x\in G_n}(x)\otimes f(x^{-1}).
\end{align*}
Then we have $gf\mapsto \sum_{x\in G_n}(gx)\otimes f(x^{-1})$
and $^gf\mapsto \sum_{x\in G_n}(g^{-1}x)\otimes gf(x^{-1})$.
If we put all
this together we obtain that the corestriction is induced 
by 
\[
\pi\otimes \pr:\Z_p[G_m]\otimes_{\Z_p}T/p^m\to \Z_p[G_n]\otimes_{\Z_p}T/p^n
\]
where $\pi:\Z_p[G_m]\to \Z_p[G_n]$ is the canonical surjection
(integration over the fibers) and $\pr :T/p^m\to T/p^n$
the canonical projection. This proves that 
\[
\prolim_nH^i(\Oh_{K_n}[1/S], T)\isom H^i(\Oh_{K}[1/S],\LG\otimes_{\Z_p}T).
\]
The $G$ action on $\LG\otimes_{\Z_p}T$ is only via the first factor
so that
\[
H^i(\Oh_{K}[1/S],\LG\otimes_{\Z_p}T)\isom 
H^i(\Oh_{K}[1/S],\LG)\otimes_{\Z_p}T.
\]
Note that the $\LG$ action is induced by the $G$ action on $T$ and is
by multiplication with the inverse on $\LG$. 
This proves the proposition.
\end{appendix}

Jens Hornbostel, NWF I - Mathematik, Universit{\"a}t Regensburg, 
93040 Regensburg, Germany, jens.hornbostel@mathematik.uni-regensburg.de 

Guido Kings,  NWF I - Mathematik, Universit{\"a}t Regensburg, 
93040 Regensburg, Germany, guido.kings@mathematik.uni-regensburg.de  

\end{document}